\documentclass[a4paper,12pt,leqno]{article}
\usepackage{latexsym}
\usepackage[all]{xy}
\usepackage{amssymb} % \mathbb
\usepackage{amsmath} % \mathbb
\usepackage{theorem}
\usepackage{colortbl} 
\usepackage{graphicx}
\def\Z{{\mathbb{Z}}}% \Z == \mathbb{Z}
\def\R{{\mathbb{R}}}% \R == \mathbb{R}
\def\A{{\mathcal{A}}}% \A == \mathcal{A}
\def\B{{\mathcal{B}}}% \B == \mathcal{B}
\def\calS{{\mathcal{S}}}% \S == \mathcal{S}
\DeclareMathOperator{\codim}{codim}

\DeclareMathOperator{\Der}{Der}

\DeclareMathOperator{\Shi}{\calS}
\DeclareMathOperator{\Cat}{Cat}
\DeclareMathOperator{\Ht}{ht}

\DeclareMathOperator{\wHt}{\widetilde{ht}}

\numberwithin{equation}{section}

\newcommand{\qed}{\hfill$\square$}

\theoremstyle{break}

\newtheorem{thrm}{Theorem}[section]
\newtheorem{prop}[thrm]{Proposition}

\newtheorem{lem}[thrm]{Lemma}
\newtheorem{defi}[thrm]{Definition}
\newtheorem{rmk}[thrm]{Remark}

\title{Free filtrations of affine Weyl arrangements\\
and the ideal-Shi arrangements}
\author{
Takuro Abe
\footnote
{
Department of Mechanical Engineering and Science,
Kyoto University,
Kyoto 606-8501, Japan.
abe.takuro.4c@kyoto-u.ac.jp
}
and Hiroaki Terao
\footnote{
Office of International Affairs,
Hokkaido University,
Sapporo 060-0815, Japan.
terao@math.sci.hokudai.ac.jp
}
}

\begin{document}

\maketitle

\begin{abstract}
In this article we prove that the ideal-Shi arrangements are
free central
arrangements of hyperplanes
satisfying the dual-partition formula.
Then it immediately follows that
there exists a saturated free filtration of the cone of
any affine Weyl 
arrangement such that 
each filter
is a free subarrangement satisfying
the dual-partition formula. 
This generalizes the main result in 
\cite{ABCHT} which affirmatively settled 
a conjecture by Sommers and Tymoczko \cite{SomTym}.
\end{abstract}

%\keywords{arrangement of hyperplanes.
%root system,
%Weyl arrangement,
%free  arrangement, 
%Shi arrangement,
%ideals
%}

%\subclass{32S22; 17B22}

\section{Introduction}
Let 
$\Phi$ be an irreducible crystallographic
root system of rank $\ell$ over the real number field
$\R$. Let $\Delta=\{\alpha_1,\ldots,\alpha_\ell\}$ be a simple system of $\Phi$ and $\Phi^+$ the corresponding positive 
system. An {\bf
ideal} $I \subseteq \Phi^+$ is a 
set such that if $\alpha \in I,\ \beta\in \Phi^+$ with 
$\alpha-\beta \in \sum_{i=1}^\ell \Z_{\ge 0} \alpha_i$, then $\beta \in I$. 
For any subset $\Sigma$ of $\Phi^{+} $,
define
$\A(\Sigma) := \{H_{\alpha} ~|~\alpha\in\Sigma \}$,
where $H_{\alpha} $ is the hyperplane perpendicular to $\alpha$
in the $\ell$-dimensional Euclidean space.

Let $V$ be the $(\ell+1)$-dimensional Euclidean space with a basis
$\Delta \cup \{z\}$,
where $z$ is a unit vector perpendicular to
each $\alpha_{i} \ (1\le i \le \ell)$. 
We usually identify $V$ with its dual space $V^{*} $ by
the inner product.
For $j\in\Z$
and
$\alpha \in \Phi^+$, define a hyperplane
$H_\alpha^j:=\{\alpha-jz=0\}$
in $V$. Let $H_{z}$ denote the hyperplane
in $V$ defined by $\{z=0\}$.

\begin{defi} 
{\rm For $k\in\Z_{> 0}$ and 
an ideal $I\subseteq \Phi^+$,
 define the
\textbf{ideal-Shi} arrangements in $V$ by 
\begin{eqnarray*}
\Shi^k_{+I}:&=&
\{
H_\alpha^{j} \mid \alpha \in \Phi^{+},  -k+1\le j\le k
\} 
\cup
\{H_{z} \}
\cup 
\{H_\alpha^{-k} \mid \alpha \in I\},
\\
\Shi^k_{-I}:&=&
\left(
\{
H_\alpha^{j} \mid \alpha \in \Phi^{+},  -k+1\le j\le k
\} \cup
\{H_{z} \}
\right)
 \setminus \ \{H_\alpha^{k} \mid \alpha \in I\}.
\end{eqnarray*}
}
\end{defi} 

The following Theorems \ref{main}, \ref{dp} and \ref{fp} 
are the main theorems of this article.
They
concern  the {\bf freeness}
 and the {\bf exponents} of the
ideal-Shi arrangements.  (See $\S$ 2 
for the terminology of the theory of (free) arrangements of hyperplanes.)

\begin{thrm}
All the ideal-Shi arrangements $\Shi^k_{\pm I}$ are free.
\label{main}
\end{thrm}

To determine the exponents of ideal-Shi arrangements, we define
the \textbf{dual partition}:
Let 
$\Sigma$ 
be a finite set of vectors in a $d$-dimensional vector space 
and $f:\Sigma \rightarrow \Z_{>0}$ be a function such that 
$f_i:=|f^{-1}(i)| \le f_{i-1}=|f^{-1}(i-1)|\leq d$ for all $ i \in \Z_{>1}$. Define
$m :=\max_{s \in \Sigma} \{f(s)\}$. 
%subspace spanned by $\Sigma$ 
Then the dual partition of the pair
$(\Sigma, f)$ is the set of integers 
$$((0)^{d-f_{1} }, (1)^{f_1-f_2},(2)^{f_2-f_3},\ldots,(m -1)^{f_{m-1}-f_m},
(m)^{f_{m}}),$$ 
where 
$(a)^b\ (a, b \in \Z_{\ge 0})$ indicates that there are $b$
copies of $a$.
The most famous example of the dual partition is the case
when $\Sigma=\Phi^{+} $ 
and
$f=\Ht : \Phi^{+} \rightarrow \Z_{>0}$, where $\Ht$ is the height function
defined by
$
\Ht
(\sum_{i=1}^{\ell} c_{i}\alpha_{i})
=
\sum_{i=1}^{\ell} c_{i}.
$
In this case,
the dual partition of 
the pair $(\Phi^{+}, \Ht)$
is equal to 
the exponents of the root system 
\cite{St} \cite{K} \cite{M}.
This remarkable dual-partition formula is generalized
to the ideal-Shi arrangements as follows:

\begin{thrm}
Let $k\in\Z_{> 0}$ and 
$I\subseteq \Phi^+$ be
an ideal. 
Denote the Coxeter number of $\Phi$ by $h$.
For $\alpha \in \Phi^+$ and $j\in\Z$, define the extended height function
$\widetilde{\Ht}$ by
$$
\widetilde{\Ht}(\alpha-jz)=
\begin{cases}
-\Ht(\alpha)+jh+1& \text{if~}j > 0,\\
\Ht(\alpha)-jh & \text{if~}j  \le 0.
\end{cases}
$$
We also define $\widetilde{\Ht}(z)=1.$ 
Then 

(1) the exponents of $\Shi^k_{+I}$ is the 
dual partition of the pair
$$
(\{\alpha-jz \mid \alpha\in\Phi^{+}, -k+1\leq j\leq k\}
\cup
\{\alpha+kz \mid \alpha \in I\}
\cup
\{z\}, \widetilde{\Ht}),$$

(2) the exponents of $\Shi^k_{-I}$ is the 
dual partition of the pair
$$
(\{\alpha-jz \mid \alpha\in\Phi^{+}, -k+1\leq j\leq k\}
\setminus
\{\alpha-kz \mid \alpha\in I\}
\cup\{z\}, 
\widetilde{\Ht}).$$
\label{dp}
\end{thrm}

\begin{rmk}
{\rm Note that
$\Shi^{k}_{\pm I}$ are equal to the cones 
(i.e., \cite[Definition 1.15]{OT})
of the
(extended and generalized) Shi arrangement
$\Shi^{k}  $ when $I$ is the empty set.
Also note that
$\Shi^{k}_{+I}$
is equal to the cones of the 
(extended and generalized) Catalan arrangement
 $ \Cat^{k}  $ when $I=\Phi^{+} $.
In these cases, Theorems~\ref{main} and \ref{dp} 
had been conjectured by
Edelman and Reiner in \cite{EdRei} before
they were proved by Yoshinaga in \cite{Y}.
}
\end{rmk}
%
%
%
%
%As an application of Theorems \ref{main}
%and \ref{dp}, 
%we can show the freeness of 
%infinitely many
% subarrangements of the affine Weyl arrangements. 	

For a
central arrangement $\A$ of
countably infinite hyperplanes,
we say that a filtration
\[
\A_{1} \subseteq 
\A_{2} \subseteq 
\dots
\subseteq
\A_{i} \subseteq 
\dots
\ \ \text{with} \ \
\A
=\bigcup_{i=1}^{\infty} \A_{i}  
\]
of $\A$ is said to be {\bf saturated}
if $|\A_{i}| = i$ for any $i\in \Z_{> 0}$.
We also say that
the filtration is {\bf free} if 
each $\A_{i}$ is a free arrangement.
Then 
Theorems \ref{main} and \ref{dp} immediately
imply the following:

\begin{thrm}
For a root system
$\Phi$, fix a linear order
$(\alpha_1,\ldots,\alpha_n)$ on the set $\Phi^{+}$ 
of positive roots
in such a way that  
$\{\alpha_i\}_{i=1}^k$ is an ideal of $\Phi^+$ for any $1 \le k \le n$. 
Define 
$$
H_{i}^j
:=H_{\alpha_i}^j
=
\{\alpha_{i} - jz = 0\}\ \ \
(j \in \Z,\ 
1\le i\le n)$$
and
\[
K_{p} :=
\begin{cases} 
H_{r}^{-q}  & \text{ if } \ 1\leq r \leq n,\\
H_{2n+1-r}^{q+1}  & \text{ if } \ n+1\leq r \leq 2n,\\
\end{cases}
\]
where $p \in \Z_{>0} $ with 
$p=2nq+r \ \ (1\leq r \leq 2n, \ q\in\Z_{\geq 0} )$. 
Let
$$
\A_{i} :=\{H_{z}, K_{1}, K_{2},\dots, K_{i-1}\}
\ \ \ (i\in \Z_{>0}). 
$$ 
Then the filtration
$
\A_{1} \subseteq 
\dots
\subseteq
\A_{i} \subseteq 
\dots
$
of
the cone 
\[
\A_{\infty} (\Phi)
:=
\{H_{z} 
\}
\cup
\{
H_i^j
\mid
j \in \Z,\ 
1\le i\le n
\}
\]
of the affine Weyl arrangement
is saturated and free.
Moreover, the exponents of 
$\A_{i}$ is the 
dual partition of the pair
$$
(\{z\}\cup\{\alpha-jz \mid \alpha\in\Phi^{+}, j\in\Z, \{\alpha-jz=0\}\in \A_{i}\}, 
\widetilde{\Ht}).$$

\label{fp} 
\end{thrm} 

%
%A typical example of the linear order $\preceq$ in Corollay~\ref{fp} 
%is a linear order of
%$\Phi^{+} = \{\alpha_{1}, \dots, \alpha_{n}\}$ satisfying
%$$
%\Ht(\alpha_{i})
%\leq
%\Ht(\alpha_{j})\
%\
%\text{whenever}
%\
%\
%i\leq j.
%$$

In Theorems \ref{main}
and \ref{dp}, we considered two types of ideal-Shi
arrangements $\Shi^k_{+I}$ and 
$\Shi^k_{-I}$. In fact, 
for an arbitrary subset $\Sigma$ of $\Phi^+$, 
we have the following theorem asserting 
a symmetry of the freeness and the exponents
with respect to $\Shi^k$. 
%, which seems to have some symmetry. In fact, the following 
%duality theorem on freeness holds.

\begin{thrm}	
For $k\in\Z_{> 0}$ and 
an arbitrary subset $\Sigma \subset \Phi^+$,
 define
\begin{eqnarray*}
\Shi^k_{+\Sigma}:&=&
\{
H_\alpha^{j} \mid \alpha \in \Phi^{+},  -k+1\le j\le k
\} 
\cup
\{H_{z} \}
\cup 
\{H_\alpha^{-k} \mid \alpha \in \Sigma\},
\\
\Shi^k_{-\Sigma}:&=&
\left(
\{
H_\alpha^{j} \mid \alpha \in \Phi^{+},  -k+1\le j\le k
\} \cup
\{H_{z} \}
\right)
 \setminus \ \{H_\alpha^{k} \mid \alpha \in \Sigma\}.
\end{eqnarray*}
Then 
the arrangement $\Shi^k_{+\Sigma}$ is free with exponents 
$
(1, kh+m_1,\ldots,kh+m_\ell)
$ if and only if 
the arrangement $\Shi^k_{-\Sigma}$ is free 
with exponents 
$
(1, kh-m_1,\ldots,kh-m_\ell).$ 
In this case, the arrangement 
$\A(\Sigma)$
 is also free with 
exponents $(m_1,\ldots,m_\ell)$.
\label{dual}
\end{thrm}

The organization of this article is as follows. In \S2 we 
present the four earlier results 
which will play important roles in the subsequent sections. 
They are
the two freeness criteria (Theorem \ref{yoshinaga} in \cite{Y2}
and
Theorem \ref{free} in 
\cite{AY2}), 
the 
ideal free theorem (Theorem \ref{ABCHT}, \cite{ABCHT}), 
and
the shift isomorphism (Theorem \ref{AY} in  \cite{AY}). 
In \S3 we prove Theorem \ref{dual} after 
we study root systems of rank two. 
In \S4 we prove Theorems \ref{main}
and
\ref{dp}.

\section{Preliminaries}
In this section
let $\A$ be a central arrangement of hyperplanes 
in $V=\R^n$, i.e., 
a finite set of hyperplanes of $V$ going through the origin. 
For each $H \in \A$ fix 
a linear form $\alpha_H \in V^*$ such that $\ker \alpha_H=H$. 
Let $Q(\A):=\prod_{H \in \A} \alpha_H$. 
An \textbf{intersection 
poset} $L(\A)$ is defined by 
\begin{eqnarray*}
L(\A):=\{\bigcap_{H \in \B} H \mid \B \subseteq \A\},\ \ \ 
L_k(\A):=\{X \in L(\A) \mid \codim X=k\}\ \ (k \in \Z_{\ge 0}).
\end{eqnarray*}
Then,
ordered by reverse inclusion,
$L(\A)$ is a poset 
with the minimum element $V$. 
The 
{\bf characteristic polynomial} 
$\chi(\A)$ is defined by 
$$
\chi(\A,t):=\sum_{X \in L(\A)} \mu(X) t^{\dim X},
$$
where the {\bf M\"obius function}
$\mu:L(\A) \rightarrow \Z$ is defined by 
\[
\mu(X)=\begin{cases}
1 & (X=V),\\ 
-\sum_{V \supset Y \supsetneq X} \mu(Y) & (X \neq V).
\end{cases}
\]
Since $\A$ is central, it is known that 
$\chi(\A,t)$ is divisible by $t-1$. Define $\chi_0(\A,t):=\chi(\A,t)/(t-1)$. 

For $X\in L(\A)$, define
\[
\A_{X} := \{H\in\A ~|~ X\subseteq H\},
\ \ \ 
\A^{X} := \{K\cap X  ~|~ K\in \A\setminus \A_{X} \}.
\]

Let $S=S(V^*)$ be a symmetric algebra of $V^*$, $\Der S$ the 
derivation module of $S$ and 
$\Omega^q_S$ the $S$-module of regular differential $q$-forms. Define 
\begin{eqnarray*}
D(\A)&=&
\{\theta \in \Der S \mid \theta(\alpha_H) \in S \alpha_H \ \text{for all~} H \in \A\},\\
\Omega^q(\A)&=&\{
\omega \in (1/Q(\A))\Omega^1_S \mid 
 Q(\A) \omega \wedge d \alpha_H \in \alpha_H \Omega^{q+1}_S\ \text{for all~}
 H \in \A\}.
\end{eqnarray*}

It is known 
(e. g., \cite{OT})
that the $S$-modules $D(\A)$ and $\Omega^1(\A)  $ 
are dual to each other.
We say that
$\A$ is \textbf{free with exponents $\exp(\A)=(d_1,\ldots,d_n)$} if
there are homogeneous derivations $\theta_1,\ldots,\theta_n \in D(\A)$ such that 
$D(\A)=\oplus_{i=1}^n S  \theta_i$ with $\deg \theta_i=d_i\ (i=1,\ldots,n)$. 
By the duality above, $\A$ is free with exponents $(d_1,\ldots,d_n)$ if and only if 
$\Omega^1(\A)$ is a free $S$-module of rank $n$ with homogeneous basis $\omega_1,\ldots,
\omega_n$ such that $\deg \omega_i=-d_i \ (i=1,\ldots,n)$.

Finally, let us introduce four
key results to prove Theorem \ref{main}. To state them, let us introduce 
multiarrangements. For a central arrangement $\A$ and $m:\A \rightarrow \Z_{>0}$, the pair 
$(\A,m)$ is called a {\bf multiarrangement}. Define 
$$
D(\A,m):=\{\theta \in \Der S \mid 
\theta(\alpha_H) \in S\alpha_H^{m(H)}\ \text{for all~} H \in \A\}.
$$
Also, the freeness of $(\A,m)$ 
and the exponents $\exp(\A,m)$ can be defined in the same way as the freeness of 
$\A$ and $\exp(\A)$. 
For a fixed $H_0 \in \A$, define $m_0:\A^{H_0} \rightarrow \Z_{>0}$ by 
$$
m_0(H \cap H_0):=
|\{K \in \A \setminus \{H_0\} \mid K \cap H_0=H \cap H_0\}|.
$$ 
The multiarrangement 
$(\A^{H_0},m_0)$ is called the \textbf{Ziegler restriction} of $\A$ 
onto $H_0$. 
If $\A$ is free with $\exp(\A) = (1, d_{2}, \dots, d_{n})$,
then  $(\A^{H_0},m_0)$ is free with 
$\exp(\A^{H_0},m_0) = (d_{2}, \dots, d_{n})$.
For $D_0(\A):=\{\theta \in D(\A) \mid \theta(\alpha_{H_0})=0\}$, 
define the 
\textbf{Ziegler restriction map} $D_0(\A) \rightarrow D(\A^{H_0},m_0)$ as the 
restriction of a
derivation onto $H_0$. For details, see \cite{Z}.

\begin{thrm}[\cite{Y2}, Theorem 3.2]
Let $\A$ be a central arrangement in $\R^3$,
$H_0 \in \A$ and $(\A'',m)$ the Ziegler restriction of $\A$ onto $H_0$. 
Let $\exp(\A'',m)=(d_1,d_2)$. Then $\A$ is free with $\exp(\A)=(1,d_1,d_2)$ if 
and only if 
$\chi_0(\A,0)=d_1d_2$.
\label{yoshinaga}
\end{thrm}

\begin{thrm}[\cite{AY2}, Theorem 4.1]
Let $\A$ be a central arrangement in $\R^n$
$(n>3)$ 
and fix $H_0 \in \A$. Let $(\A'',m)$ be the 
Ziegler restriction of $\A$ onto $H_0$. Assume that 

(1)
$(\A'',m)$ is free, and 

(2)
$\A_X$ is free for any $X \in L_3(\A)$ with $X \subset H_0$.

\noindent
Then $\A$ is free.
\label{free}
\end{thrm}

We use the notation from \S 1:
let $\Phi$ be an irreducible crystallographic root system
of rank $\ell$.

\begin{thrm}[Shift isomorphism, \cite{AY}, Corollary 12]
Let  $k\in\Z_{>0} $,
$\A:=\A(\Phi^{+} )$
and $m: \A \rightarrow \{0,1\}$ be a multiplicity. 
Then there exist isomorphisms of $S$-modules 
\begin{eqnarray*}
D(\A,m) \rightarrow D(\A,2k+m),\ \ \ \
\Omega^1(\A,m) \rightarrow D(\A,2k-m).
\end{eqnarray*}
Hence if $(\A,m)$ is free with $\exp(\A,m)=(m_1,\ldots,m_{\ell})$, then 
$(\A,2k\pm m)$ is also free with $\exp(\A,2k \pm m)=
%((kh)^{\ell}) 
%\pm (m_1,\ldots,m_{\ell})=
(kh\pm m_1,\ldots,kh\pm m_{\ell})
$. 
\label{AY}
\end{thrm}

%
%If we define 
%then
%still hold true.
%It had been
%conjectured
%by 
In 
\cite{SomTym}
Sommers and Tymoczko posed the
conjecture corresponding to
Theorems \ref{main} and \ref{dp} 
for $k=0$.  The conjecture
was affirmatively settled
as follows:

\begin{thrm}[Ideal free theorem, \cite{ABCHT}, Theorem 1.1]
Let
$I \subseteq \Phi^+$ be an ideal.
Then 
$\A(I)=\{H_\alpha  \mid \alpha \in I\}$
is a free arrangement and its exponents 
$(m_{1}(I), m_{2}(I), \dots, m_{\ell}(I))$ 
are equal to 
the dual partition of the pair $(I, \Ht)$,
%where $\Ht$ denotes the height function of positive roots.
where $\Ht$ is the height function of positive roots.
\label{ABCHT}
\end{thrm}

\section{Proof of Theorem \ref{dual}}

In this section 
we continue to
use the notation from \S 1.
Before the proof of 
Theorem \ref{dual},
we will verify
Lemma \ref{key0} and then
prove 
Proposition \ref{key1}
 which is a key to the proof of 
Theorem \ref{dual}.
In Lemma \ref{key0} and Proposition \ref{key1}, let 
$\Phi$ denote an irreducible crystallographic root system of rank two
(i.e., $\Phi = A_{2}, B_{2}  $ or $G_{2} $).
Recall that $\Delta$ is a simple system of $\Phi$.
 Fix $k\in\Z_{>0}$.

\begin{lem}
For $\alpha, \beta\in\Phi^{+} \ (\alpha\neq \beta)$,
let $p_{-} := H_{\alpha}^{k}\cap H_{\beta}^{k}$ .

(1) 
If 
$\Delta = \{\alpha, \beta\}$,
then
$\{H_{\gamma}^{s}~|~ 
\gamma\in\Phi^{+},
-k\le s \le k,
p_{-} \subset H_{\gamma}^{s}
\}
=\{H_{\alpha}^{k}, H_{\beta}^{k}\}$,

(2) 
if $\Delta \neq \{\alpha, \beta\}$,
then there exists
$\gamma\in\Phi^{+} $ 
such that
$p_{-} \subset H_{\gamma}^{0}$.

\noindent
These two results hold true also for
$p_{+} := H_{\alpha}^{-k}\cap H_{\beta}^{-k}$.  
\label{key0}
\end{lem}

\noindent
\textbf{Proof}. 
(1) Assume that 
$\Delta = \{\alpha, \beta\}$
and that
$p_{-} \subset H_{\gamma}^{s}  $ for some
$\gamma\in \Phi^{+} $ and some $s$ with
$-k\le s \le k$.
Then we have
\[
a(\alpha-kz)+b(\beta-kz)=\gamma-sz
\]
 for some nonzero
rational numbers $a, b.$ 
Since $a\alpha+b\beta=\gamma$, one has
$\{a, b\}\subset\Z_{>0} $.  Thus 
$s=ak+bk=(a+b)k>k$,   
which is a contradiction.

(2) If $\Delta \neq \{\alpha, \beta\}$, by case-by-case arguments 
for $A_{2}, B_{2}, G_{2}$, we have
$\alpha-\beta \in \Z \gamma$ for some $\gamma\in\Phi^{+} $.  
This implies $p_{-} = H_{\alpha}^{k}\cap H_{\beta}^{k}\subset H_{\gamma}^{0}$
because $(\alpha-kz)-(\beta-kz)\in \Z \gamma$.  

In the case of $p_{+}$, the parallel proof works. 
\qed

\medskip

For an arbitrary arrangement $\A$ and a hyperplane $H_{0} $, define
\[
\A\cap H_{0} := 
\{
K\cap H_{0} ~|~ K\in\A, K\neq H_{0} 
\}.
\]
Then $\A\cap H_{0}$  is an arrangement in $H_{0}. $ 

\begin{prop}
Let
$\Sigma\subseteq \Phi^{+}$
and
 $\alpha \in \Phi^{+}\setminus\Sigma$. 
Then

(1)
\[
\left|
\Shi^{k}_{+\Sigma} \cap H_{\alpha}^{-k}   
\right|
=
\begin{cases} 
kh+1 \ \ \ \text{if~}\alpha\in\Delta, \Sigma\cap \Delta =\emptyset,\\
kh+2 \ \ \ \text{otherwise}, 
\end{cases} 
\]

(2)
\[
\left|
\Shi^{k}_{-\Sigma} \cap H_{\alpha}^{k}   
\right|
=
\begin{cases} 
kh+1 \ \ \ \text{if~}\alpha\in\Delta, \Sigma\cap \Delta =\emptyset,\\
kh~~~~~ \ \ \ \text{otherwise}. 
\end{cases} 
\]
\label{key1} 
\end{prop}

\noindent
\textbf{Proof}. 
When $\Sigma=\emptyset$, 
by directly counting the intersections, 
we get the following equalities (\cite{AT}): 
\begin{equation}
\label{count1} 
\left|
\Shi^{k} \cap H_{\alpha}^{-k}   
\right|
=
\begin{cases} 
kh+1 \ \ \ \text{if~}\alpha\in\Delta,\\
kh+2 \ \ \ \text{otherwise}, 
\end{cases} 
\end{equation} 
\[
\left|
\Shi^{k} \cap H_{\alpha}^{k}   
\right|
=
\begin{cases} 
kh+1 \ \ \ \text{if~}\alpha\in\Delta,\\
kh~~~~~ \ \ \ \text{otherwise}. 
\end{cases} 
\]

(1)
Consider 
the difference set
\[
D_{+} 
:=
\left(
\Shi^{k}_{+\Sigma}
\cap H_{\alpha}^{-k}   
\right)
\setminus
\left(
\Shi^{k} 
\cap H_{\alpha}^{-k}   
\right).
\]
{\it
Case 1.} Suppose $\alpha\not\in \Delta$. 
Let $\beta\in \Sigma.$ 
Then 
$
p_{+} =
H_{\beta}^{-k} 
\cap
H_{\alpha}^{-k} 
\subset
H_{\gamma}^{0}  
$
for some $\gamma\in\Phi^{+} $ 
by Lemma \ref{key0} (2).  This implies
\[
p_{+} =
H_{\beta}^{-k} 
\cap
H_{\alpha}^{-k} 
=
H_{\gamma}^{0}  
\cap
H_{\alpha}^{-k} 
\in
\Shi^{k} 
\cap
H_{\alpha}^{-k}.  
\]
Therefore 
$D_{+} = \emptyset$.

{\it
Case 2.} 
Suppose $\alpha\in \Delta$
and $\Sigma\cap\Delta=\emptyset$.
Then an arbitrary root $\beta\in\Sigma$ is non-simple.
Thus 
$
p_{+} =
H_{\beta}^{-k} 
\cap
H_{\alpha}^{-k} 
\subset
H_{\gamma}^{0}  
$
for some $\gamma\in\Phi^{+} $ 
by Lemma \ref{key0} (2).  This implies
\[
p_{+} =
H_{\beta}^{-k} 
\cap
H_{\alpha}^{-k} 
=
H_{\gamma}^{0}  
\cap
H_{\alpha}^{-k} 
\in
\Shi^{k} 
\cap
H_{\alpha}^{-k}.  
\]
Therefore 
$D_{+} = \emptyset$.

{\it
Case 3.} 
Suppose $\alpha\in \Delta$ and
$\Sigma\cap\Delta \neq\emptyset$. 
Then we may express $\Delta
=
\{\alpha, \beta\}$ 
with $\beta\in\Sigma$.
By Lemma \ref{key0} (1)
we have
\[
p_{+} =
H_{\beta}^{-k} 
\cap
H_{\alpha}^{-k} 
\in
\Shi^{k}_{+\Sigma}  
\cap
H_{\alpha}^{-k},
\ \ \ \
p_{+} =
H_{\beta}^{-k} 
\cap
H_{\alpha}^{-k} 
\not\in
\Shi^{k} 
\cap
H_{\alpha}^{-k}.  
\]
Therefore 
$D_{+} = \{p_{+} \}$.

Combining (\ref{count1}) with the three cases above, we get
\[
\left|
\Shi^{k}_{+\Sigma} \cap H_{\alpha}^{-k}   
\right|
=
\begin{cases} 
\left|
\Shi^{k} \cap H_{\alpha}^{-k}   
\right|
=
kh+1~~~~
 \ \ \ \text{if~}\alpha\in\Delta, \Sigma\cap \Delta =\emptyset,\\
\left|
\Shi^{k} \cap H_{\alpha}^{-k}   
\right|+1
=
kh+2~ \ \text{if~}\alpha\in\Delta, \Sigma\cap \Delta \neq\emptyset,\\
\left|
\Shi^{k} \cap H_{\alpha}^{-k}   
\right|
=
kh+2~~~~ \ \ \ \text{otherwise}. 
\end{cases} 
\]
This proves (1).

As for (2), the parallel proof works
if we use
the difference set
\[
D_{-} 
:=
\left(
\Shi^{k} 
\cap H_{\alpha}^{k}   
\right)
\setminus
\left(
\Shi_{-\Sigma}^{k} 
\cap H_{\alpha}^{k}   
\right)
\]
instead of $D_{+}. $ 
\qed

\medskip

\noindent
{\bf Proof of Theorem \ref{dual}.}

Recall that the Shi arrangements $\Shi^{k}$ are free with
%$\exp(\Shi^{k}) = (1, (kh)^{\ell})$ and
$\exp(\Shi^{k}) = (1, kh, \dots, kh)$ and
$\chi(\Shi^{k}, t) = (t-1)(t-kh)^{\ell}$
by \cite{Y}. 

\smallskip

{\it Claim.}
Assume $\ell=2$.  
Then 
$\Shi^k_{\pm \Sigma}$ is free if and only if
either $\Sigma=\emptyset$ or
$\Sigma\cap\Delta\neq\emptyset$. 
Therefore Theorem \ref{dual} holds true when $\ell=2.$ 

\smallskip

Let us verify \textit{Claim}.  Recall $\chi_{0} (\A, t) := \chi(\A, t)/(t-1)$ from \S 2. 
We will apply Theorem \ref{yoshinaga}.
Define $\zeta(\Sigma):=\chi_{0} (\Shi_{+\Sigma}^{k}, 0)$ for an
arbitrary subset $\Sigma$ of $\Phi_{+} $. 
Note that $\zeta(\emptyset)=(kh)^{2} $. 
For $\alpha\in\Phi^{+}\setminus \Sigma$ one has
\begin{equation}
\zeta (\Sigma\cup\{\alpha\})
=
\zeta (\Sigma)
-
\chi_{0} (\Shi_{+\Sigma}\cap H_{\alpha}^{-k}, 0)
\label{delres} 
\end{equation}
because of the deletion-restriction formula for $\chi$ 
(i. e., \cite[Corollary 2.57]{OT}). 
Since $\Shi_{+\Sigma}\cap H_{\alpha}^{-k}$ is an arrangement
in the real $2$-dimensional space $H_{\alpha}^{-k}$,
we obtain
$$\chi_{0} (\Shi_{+\Sigma}\cap H_{\alpha}^{-k}, 0)
=1-
\left|
\Shi_{+\Sigma}\cap H_{\alpha}^{-k}\right|.$$   
Thanks to (\ref{count1}) and (\ref{delres}), for $\alpha\in\Delta$, we have 
\[
\zeta(\{\alpha\})=
\begin{cases} 
(kh)^{2} + kh
 \ \ \ \ \ \ \ \   \ \ \  \text{if~}\alpha\in\Delta,\\
(kh)^{2} + (kh+1)
\ \ \ \text{otherwise}. 
\end{cases} 
\]
Similarly we may verify
\begin{equation} 
\zeta(\Sigma)=
\begin{cases} 
(kh)^{2} + kh + (kh+1)(|\Sigma|-1) 
%= (kh+1)(kh+|\Sigma|-1)
 \ \ \ \text{if~}\Sigma\cap\Delta\neq\emptyset,\\
(kh)^{2} + (kh+1)|\Sigma|
\ \ \ \ \ \ \ \ \ \ \ \ \ \ \ \ \ \ \text{otherwise}
\label{zeta} 
\end{cases} 
\end{equation}
by applying Proposition \ref{key1} 
and (\ref{delres})
repeatedly.

Now we will apply Theorem \ref{AY}. 
Let $\A=\A(\Phi^{+} )$.   Suppose that
$m: = {\bf 1}_{\Sigma}$
is the indicator function of 
$\A(\Sigma)$ in $\A$. 
Note that the Ziegler restriction of $\Shi_{+\Sigma}^{k}$ onto
 $H_{z} $ is $(\A, 2k+m)$. 
Let
$(d_{1}, d_{2})
:=
\exp
(\A, 2k+m)$.
Define $\zeta'(\Sigma):=d_{1} d_{2} $. 
Note that
$\zeta'(\emptyset) = (kh)^{2} $. 
Suppose $\Sigma\neq\emptyset$.
 Since
$
\exp(\A, m)=
(1, |\Sigma|-1),
$
Theorem \ref{AY} gives
\begin{equation*} 
\exp(\A, 2k+m)
=
\begin{cases} 
(kh+1, kh+|\Sigma|-1)
 \ \ \ \ &\text{if~}\Sigma\neq\emptyset,\\
(kh, kh)
&\text{if~}\Sigma=\emptyset.
\end{cases} 
\end{equation*} 
Thus we obtain
\begin{equation} 
\zeta'(\Sigma)
=
\begin{cases} 
(kh+1)(kh+|\Sigma|-1)
 \ \ \ \ &\text{if~}\Sigma\neq\emptyset,\\
(kh)^{2}
&\text{if~}\Sigma=\emptyset.
\end{cases} 
\label{zetadash} 
\end{equation} 
Comparing the equations (\ref{zeta}) and (\ref{zetadash}), we may
conclude that $\zeta(\Sigma)=\zeta'(\Sigma)$ if and only if
either $\Sigma=\emptyset$ or
$\Sigma\cap\Delta\neq\emptyset.$  
This shows {\em Claim} for $\Shi_{+\Sigma}^{k}$
by Theorem \ref{yoshinaga}.
It is not hard to see that
the parallel proof works for $\Shi_{-\Sigma}^{k}$.

%
%Then Proposition \ref{key1} and the addition/deletion theorem in \cite{T} 
%show that $\Shi^k_{\pm \Sigma}$ is free.
%% if $\Sigma$ contains a simple root. 
%Assume that $\Sigma$ does not contain any simple root. Then 
%Lemma \ref{key0} shows that 
%$\chi_0(\Shi^k_{\pm \Sigma},0)
%=(kh)^2 \pm (|\Sigma|kh +|\Sigma|)>(kh \pm 1)(kh \pm (|\Sigma|-1))$. 
%On the other hand, the exponents of the Ziegler restriction of 
%$\Shi^k_{\pm \Sigma}$ onto $H_{z}$ are $(kh \pm 1, kh \pm (|\Sigma|-1))$ by 
%Theorem \ref{AY}. 
%Hence Theorem \ref{yoshinaga} shows that $\Shi^k_{\pm \Sigma}$ is not free, 
%which 
%completes the proof when $\ell=2$.

\smallskip

Next assume that $\ell \geq 3$.
We will apply
Theorem \ref{free}.
We still use the notation
$\A=\A(\Phi^{+} )$
and
$m = {\bf 1}_{\Sigma}$.
Then the Ziegler restriction 
of 
$\Shi^k_{\pm\Sigma}$ onto 
$H_{z} $ is equal to
$(\A, 2k\pm m)$. 
Theorem \ref{AY} 
shows that 
$(\A,2k+m)$ is free if and only if
$(\A,2k-m)$ is free. 
Let $X \in L_3(\Shi^k_{-\Sigma})$
 with $X \subset H_{z}$. It is known that $X=Y \cap H_{z} $ for some 
$Y \in L_2(\A)$ (see \cite{AT} for example). 
Note that $\Psi:=\Phi \cap Y^{\perp}$ 
is a (not necessarily irreducible)
root system of rank two.
Then $\Psi^{+} :=\Phi^{+}  \cap Y^{\perp}$ 
is a positive system of $\Psi$. 
Define 
$\B_{\pm}$ to be the restriction of 
$(\Shi^k_{\pm \Sigma})_X$
to the 2-dimensional vector space $Y^{\perp}$.
Then
$\B_{\pm} $ is equal to 
$\Shi_{\pm(\Sigma\cap\Psi^{+})}^{k}  $ when 
the entire root system is equal to $\Psi$.  
 {\em Claim} shows that
 $\B_{+}$ is free if and only if 
$\B_{-}$ is free.
(We may easily check both $\B_{+}$ 
and
$\B_{-}$ are free for $\Psi=A_{1} \times A_{1} $.)
Note that
$
(\Shi^k_{\pm \Sigma})_X
=
\B_{\pm} \times \{X\}$, where $\{X\}$ 
is a singleton arrangement in $Y$.
Therefore 
$\B_{\pm}$
is free if and only if
$(\Shi^k_{\pm\Sigma})_X$ is free. 
Hence 
%\textit{Claim} implies
 $(\Shi^k_{+\Sigma})_X$ is free if and only if 
$(\Shi^k_{-\Sigma})_X$ is free.
Now we may apply Theorem \ref{free} to conclude that
the freeness of $\Shi^k_{+\Sigma}$ is equivalent to
the freeness of $\Shi^k_{-\Sigma}$.
If they are free, then
$\A(\Sigma)=(\A,m)$ is also free
by  Theorem \ref{AY}. Let 
$\exp(\A(\Sigma)):=(m_1,\ldots,m_\ell)$. 
Then
$\exp(\A,2k\pm m)=(kh\pm m_1,\ldots,kh\pm m_\ell)$.
Since $(\A,2k\pm m)$ is the Ziegler restriction of
$\Shi^k_{\pm\Sigma}$ onto $H_{z} $, we conclude that 
$\exp(\Shi^k_{\pm \Sigma})
=(1,kh\pm m_1,\ldots, kh \pm m_\ell)$, 
which completes the proof. \qed

\section{Proofs of Theorems \ref{main} and \ref{dp} }
In this section let us prove Theorem \ref{main}. 
For that purpose, we first introduce the following 
lemma:

\begin{lem}
Let $I \subset \Phi^+$ be an ideal and 
$X \in L_3(\Shi^k_{\pm I})$ such that 
$X \subset H_{z} $. 
Let $\A=\A(\Phi^+)$. 
Choose $Y \in L_2(\A)$ such that 
$X=Y \cap H_{z} $. 
Let $\Psi:=\Phi \cap Y^{\perp}$ 
and $\Psi^{+} :=\Phi^{+}  \cap Y^{\perp}$.
Then $J:=I \cap \Psi^+$ is also an ideal of $\Psi^+$.
%
%Also, $(\Shi^k_{- I})_X=
%\{H_{z}\} \cup
%\{
%H_\alpha^j \mid -k+1 \le j \le k,\ 
%\alpha \in \Psi^+\} \setminus \{H_\alpha^k \mid \alpha \in J\}$, and 
%$(\Shi^k_{+ I})_X=
%\{H_{z}\} \cup
%\{
%H_\alpha^j \mid -k+1 \le j \le k,\ 
%\alpha \in \Psi^+\} \cup \{H_\alpha^{-k} \mid \alpha \in J\}$.
\label{ideal}
\end{lem}

\noindent
\textbf{Proof}. 
%
%As sets, the equalities on $\Shi^k_{\pm I}$ are 
%immediate. Hence it suffices to show that $J$ is an ideal in $\Psi^+_X$. 
Let $\alpha \in J$ and $\beta \in \Psi^+$ such that 
$\alpha-\beta \in \Z_{\ge 0} \gamma_1+\Z_{\ge 0} \gamma_2$, where $\{\gamma_1,\gamma_2\}$ 
is the simple system of $\Psi^+$. Since $\gamma_1$ and $\gamma_2$ are 
positive roots in $\Phi$, $\alpha \ge \beta$ in $\Phi^+$. Hence $\beta \in I$, which implies that 
$\beta \in J$. \qed

\medskip

\noindent
\textbf{Proof of Theorem \ref{main}}. 
By Theorem \ref{dual}, it suffices to show that
$\Shi^k_{+I}$ is free. 
Assume $\ell = 2$.
Note that $I\cap\Delta\neq\emptyset$ unless $I = \emptyset$. 
Hence {\em Claim} in the proof of Theorem
\ref{dual} completes the proof. 

Assume that $\ell \ge 3$.  
We apply Theorem \ref{free} to prove Theorem \ref{main}. 
For that purpose, let us verify the two conditions in 
Theorem \ref{free}.

We use the notation
$\A=\A(\Phi^{+} )$
and
$m = {\bf 1}_{I}$
as in the proof of 
Theorem \ref{dual}.
Then the Ziegler restriction 
of 
$\Shi^k_{+I}$ onto 
$H_{z} $ is equal to
$(\A, 2k + m)$. 
By 
Theorem \ref{ABCHT} 
$(\A,m)$ is free, and $D(\A, m) \simeq D(\A,2k+m)$ by 
Theorem \ref{AY}.
This verifies the condition (1) in 
Theorem \ref{free}.

Next we verify the condition (2).
Let
$X \in L_3(\Shi^k_{+ I})$ such that 
$X \subset H_{z} $. 
Choose $Y$, $\Psi$ and $\Psi^{+}$ as in
Lemma \ref{ideal}.
Define 
$\B$ to be the restriction of 
$(\Shi^k_{+I})_X$
to the 2-dimensional vector space $Y^{\perp}$.
Then
$\B$ is equal to 
$\Shi_{+(I\cap\Psi^{+})}^{k}  $ when 
the entire root system is equal to $\Psi$.  
Lemma \ref{ideal} shows that
 $I\cap\Psi^{+}\subset \Psi^{+}$ is an ideal.
This verifies the freeness of 
$
\B
$
because 
Theorem \ref{main} has been already proved when $\ell=2$. 
Recall that
$
(\Shi^k_{+I})_X
$
is free if and only if
$\B$
is free
as we saw in the proof of Theorem \ref{dual}. 
This verifies the condition (2) in 
Theorem \ref{free}.
\qed
\medskip

\noindent
\textbf{Proof of Theorem  \ref{dp}.} 

(1)
Let $k\in \Z_{>0} $. 
Define
$$A:=
\{\alpha-jz \mid \alpha\in\Phi^{+}, -k+1\leq j\leq k\},
\ \ \ \
B:=
\{\alpha+kz \mid \alpha \in I\}.
$$
Recall $
\wHt
:
A\cup
B\cup
\{z\}
\rightarrow\Z_{>0}$
from Theorem \ref{dp}.

\smallskip

\noindent
{\it Claim 1.}
$\wHt(A)\subset \left[1, kh\right]$
and
$\wHt(B)\subset \left[kh+1, (k+1)h-1\right]$.

\smallskip
Let $\alpha\in\Phi^{+} $. 
For $0<j\leq k$,
we have $1 < -\Ht(\alpha)+jh+1\leq kh$
because $1\leq \Ht(\alpha) < h$.
For $1-k\leq j\leq 0$,
we have $1\leq \Ht(\alpha)-jh < kh$
because $1\leq \Ht(\alpha)< h$.
This verifies 
$\wHt(A)\subset \left[1, kh\right]$.
Let $\beta\in I$. Similarly we may easily verify
$\wHt(B)\subset \left[kh+1, (k+1)h-1\right]$.

Consider the standard height function 
$\Ht: \Phi^{+} \longrightarrow \Z_{>0}$.
Define  
$g_{i}:= |\Ht^{-1} (i)|$.
Then $g_{h} = 0$ and $g_{1} = \ell$. 

\smallskip
\noindent
{\it Claim 2.}
$g_{i}+g_{h-i+1} = \ell $ for $1\leq i\leq h$.
\smallskip

Recall from \cite{St} \cite{K} \cite{M}
that the dual partition of the pair
$(\Phi, \Ht)$ is equal to $\exp(\A(\Phi^{+}))$:
\[
\exp(\A(\Phi^{+} ))=
((1)^{g_{1} - g_{2} }, (2)^{g_{2} - g_{3} }, 
\dots, (h-1)^{g_{h-1}} ).
\]
By the duality of the exponents, we have
$g_{i} - g_{i+1} = g_{h-i} - g_{h-i+1}$ 
and thus
$g_{i} + g_{h-i+1} = g_{i+1} + g_{h-i}$ 
for $1\leq i < h$. This implies that
the value of
$g_{i} + g_{h-i+1}$ does not depend upon 
$i$ with $1\leq i \leq h$.  It is equal to
$g_{1} +g_{h} = \ell.$ 
This verifies {\it Claim 2}.

 Note that $\wHt(z)=1$.  
For $i\in\Z_{>0}$ define $f_{i}:=|\wHt^{-1} (i)|$. 
By {\it Claim 1},  we have $f_{i} =0 $ if $(k+1)h\le i$. 

Let $1\leq i\leq kh$.
We may uniquely express $i=qh+r$ with 
$0\leq q\leq k-1$ and
$1\leq r\leq h$. 
Suppose that 
$\alpha\in\Phi^{+}$
and
$1-k\leq j\leq 0$.
Then it is not hard to see that
%, unless $r=h$, 
\[
\wHt(\alpha-jz)= i \Longleftrightarrow
j=-q \text{~and~}
\Ht(\alpha)=r.
\]
Next we assume that 
$0<j\leq k$.
Then it is not hard to see that
%, unless $r=1$, 
\[
\wHt(\alpha-jz)= i \Longleftrightarrow
j=q+1 \text{~and~}
\Ht(\alpha)=h-r+1.
\]
Now
we may conclude
\[
f_{i} =
\begin{cases} 
1+g_{1}+g_{h}=\ell+1 &\text{if~}i=1,\\
g_{r}+g_{h-r+1}=\ell &\text{if~}1 < i \leq kh.
\end{cases} 
\]
thanks to {\it Claim 2}. Thus we obtain
\[
f_{1} -f_{2}=1,  \ \ f_{2} -f_{3}=f_{3}-f_{4}=\dots=f_{kh-1}-f_{kh}=0.     
\]

Next let $kh < i< (k+1)h$.
We may uniquely express $i=kh+r$ with 
$1\leq r \leq h$. 
Suppose that 
$\beta\in I$.
Then it is not hard to see that
\[
\wHt(\beta+kz)= i \Longleftrightarrow
\Ht(\beta)=r.
\]
Define $p_{i} = |\{\beta\in I~|~\Ht(\beta)=i\}|$ 
for $1\leq i \leq h$. 
 Then we conclude that
%\[
$
f_{i} = p_{r}
$ 
when $kh < i=kh+r < (k+1)h.$
%\]
Thus we obtain
\begin{align*}
&f_{kh} -f_{kh+1}=\ell-p_{1}, 
f_{kh+2} -f_{kh+3}=p_{2} -p_{3},\\
 &f_{kh+3}-f_{kh+4}=p_{3}-p_{4},\dots, f_{kh+h-1}-f_{kh+h}=p_{h-1}-p_{h}.
\end{align*} 
Recall Theorem \ref{ABCHT} which asserts that
$$
\exp(\A(I))
=
(m_{1}(I), 
\dots,   
m_{\ell}(I)
)
=
((0)^{\ell-p_{1} },  
(1)^{p_{1} - p_{2} }, (2)^{p_{2} - p_{3} }, 
\dots).
$$
Therefore the dual partition of the pair
$(A\cup B\cup \{z\}, \wHt)$ is equal to
\[((0)^{\ell+1-f_{1} }, 
(1)^{f_{1} - f_{2} }, (2)^{f_{2} - f_{3} }, 
\dots, \dots )
=
(
1, 
kh+m_{1}(I), 
\dots,   
kh+m_{\ell}(I)
).
\]
This proves Theorem \ref{dp} (1) because of Theorem \ref{dual}.

(2)
The parallel proof works for $\Shi_{-I}^{k}$.
\qed
\medskip

{\bf Acknowledgements.}
~~
The first author is partially
supported by JSPS Grants-in-Aid for Young Scientists
(B)
No. 24740012.
The second author is partially supported by 
JSPS Grants-in-Aid for Scientific Research (A) 
No.\ 24244001.

\end{document}